\documentclass[a4paper,12pt]{article}
\usepackage{color}
\usepackage[english]{babel}
\usepackage{amsfonts}
\usepackage{amsmath}
\usepackage{amsthm}
\usepackage{amssymb}
\usepackage{bbm}
\usepackage{mathrsfs}
\usepackage{esint}
\usepackage{faktor}
\usepackage[utf8]{inputenc}
\usepackage{afterpage}
\usepackage{graphicx}

\newcommand{\N}{\mathbb{N}}

\newcommand{\R}{\mathbb{R}}

\renewcommand{\epsilon}{\varepsilon}
\renewcommand{\phi}{\varphi}

\newtheorem{lemma}{Lemma}[section]
\newtheorem{thm}[lemma]{Theorem}

\theoremstyle{definition}

\newtheorem{rmk}[lemma]{Remark}
\newtheorem{ex}[lemma]{Example}

\numberwithin{equation}{section}
\DeclareMathOperator*{\esssup}{ess \, sup}
\DeclareMathOperator*{\essinf}{ess \, inf}

\begin{document}
\title{\textbf{Some recent results on singular $p$-Laplacian systems}}
\author{
\bf Umberto Guarnotta\thanks{Corresponding author.}, Roberto Livrea\\
\small{Dipartimento di Matematica e Informatica, Universit\`a di Palermo,}\\
\small{Via Archirafi 34, 90123 Palermo, Italy}\\
\small{\it E-mail: umberto.guarnotta@unipa.it, roberto.livrea@unipa.it}\\
\mbox{}\\
\bf Salvatore A. Marano\\
\small{Dipartimento di Matematica e Informatica, Universit\`a di Catania,}\\
\small{Viale A. Doria 6, 95125 Catania, Italy}\\
\small{\it E-mail: marano@dmi.unict.it}\\
\mbox{}\\
}
\date{}
\maketitle
\begin{abstract}
Some recent existence, multiplicity, and uniqueness results for singular $p$-Laplacian systems either in bounded domains or in the whole space are presented, with a special attention to the case of convective reactions. A extensive bibliography is also provided.
\end{abstract}
\vspace{2ex}
\vspace{2ex}


\section{Introduction}
This survey paper can be divided into two parts. The first (cf. Section \ref{bound}) treats \textit{singular quasi-linear Dirichlet systems in bounded domains.} So, we study problems of the type 
\begin{equation}\label{sample1}
\left\{
\begin{alignedat}{2}
-\Delta_p u =f(x,u,v,\nabla u,\nabla v)\;\; &\mbox{in $\Omega$,}\;\; u>0\;\;
&&\mbox{in $\Omega$,}\;\; u\lfloor_{\partial\Omega}=0,\\
-\Delta_q v =g(x,u,v,\nabla u,\nabla v)\;\; &\mbox{in $\Omega$,}\;\; v>0\;\;
&&\mbox{in $\Omega$,}\;\; v\lfloor_{\partial\Omega}=0,
\end{alignedat}
\right.
\end{equation}
where $1<p,q<\infty$, the symbol $\Delta_r$ denotes the $r$-Laplace operator, namely 
\begin{equation*}
\Delta_r u:={\rm div}(|\nabla u|^{r-2}\nabla u), 
\end{equation*}
$\Omega$ is a bounded domain in $\R^N$, $N\geq 3$, having a smooth boundary $\partial\Omega$, while $f,g\in C^0(\Omega\times(\R^+)^2 \times(\R^N)^2)$ turn out singular at zero with respect the solution $(u,v)$ or even its gradient $(\nabla u,\nabla v)$.

If $p=q=2$ then various special, often non-convective, cases of \eqref{sample1} have been thoroughly investigated; see Section \ref{bound.1} below. In particular, the monograph \cite{GR1} gives a nice introduction to so-called singular Gierer-Meinhardt systems. Here, we simply make a short account on some recent existence, multiplicity, or uniqueness results when $(p,q)\neq(2,2)$, as well as the relevant technical approaches. Regarding this, let us also point out Chapter 7 of \cite{Mo}.

The second part (cf. Section \ref{wholespace}) carries out a similar analysis for \textit{singular quasi-linear systems in the whole space}. Hence, it deals with situations like
\begin{equation}\label{sample2}
\left\{
\begin{alignedat}{2}
-\Delta_p u =f(x,u,v,\nabla u,\nabla v)\;\; &\mbox{in $\R^N$,}\;\; u>0\;\;
&&\mbox{in $\R^N$,}\\
-\Delta_q v =g(x,u,v,\nabla u,\nabla v)\;\; &\mbox{in $\R^N$,}\;\; v>0\;\;
&&\mbox{in $\R^N$.}
\end{alignedat}
\right.
\end{equation}
To the best of our knowledge, no previous book or survey on problem \eqref{sample2} is already available.

Both parts exhibit four sub-sections. The first represents a historical sketch of the semi-linear setting. The next two address existence, multiplicity, and uniqueness in the non-convective case. The fourth is devoted to singular systems with convection. 

It may be useful to emphasize that \cite{Go} treats non-local singular systems, which, for the sake of brevity, are not examined here. 

We apologize in advance for possibly forgetting significant works, but the literature on \eqref{sample1}--\eqref{sample2} is by now quite extensive and our knowledge somewhat limited.
\section{Basic notation}
Let $X(\Omega)$ be a real-valued function space on a nonempty measurable set $\Omega\subseteq\R^N$. If $u_1,u_2\in X(\Omega)$ and $u_1(x)<u_2(x)$ a.e. in $\Omega$ then we simply write $u_1<u_2$. The meaning of $u_1\leq u_2$, etc. is analogous. Put 
\begin{equation*}
X(\Omega)_+:=\left\{u\in X(\Omega): u\geq 0\right\}.   
\end{equation*}
The symbol $u\in X_{\rm loc}(\Omega)$ means that $u:\Omega\to\R$ and $u\lfloor_K\in X(K)$ for all nonempty compact subset $K$ of $\Omega$. Given $1<r<+\infty$, define
$$r':=\frac{r}{r-1}\, .$$
We denote by $\lambda_{1,r}$ the first eigenvalue of the operator $-\Delta_r$ in $W^{1,r}_0(\Omega)$. If $r<N$ then
$$r^*:=\frac{Nr}{N-r}\, .$$
Let us next recall the notion and some relevant properties of the so-called Beppo Levi spaces $\mathcal{D}^{1,r}_0(\R^N)$, systematically studied for the first time by Deny and Lions \cite{DL}. Set
\begin{equation*}
\mathcal{D}^{1,r} := \left\{ z \in L^1_{\rm loc}(\R^N): |\nabla z| \in L^r(\R^N) \right\}
\end{equation*}
and write $\mathcal{R}$ for the equivalence relation that identifies two elements in $\mathcal{D}^{1,r}$ whose difference is a constant. The quotient set $\mathcal{\dot D}^{1,r}$, endowed with the norm
$$\|u\|_{1,r} :=\left(\int_{\R^N}|\nabla u(x)|^r{\rm d}x\right)^{1/r},$$
turns out complete. Indicate with $\mathcal{D}^{1,r}_0(\R^N)$ the subspace of
$\mathcal{\dot D}^{1,r}$ defined as the closure of $C^\infty_0(\R^N)$ under
$\|\cdot\|_{1,r}$, namely
$$\mathcal{D}^{1,r}_0(\R^N) := \overline{C^\infty_0(\R^N)}^{\|\cdot\|_{1,r}}.$$
$\mathcal{D}^{1,r}_0(\R^N)$, usually called Beppo Levi space,  is reflexive and continuously embeds in $L^{r^*}(\R^N)$, i.e., 
\begin{equation}\label{embedding}
\mathcal{D}^{1,r}_0(\R^N) \hookrightarrow L^{r^*}(\R^N).
\end{equation}
Consequently, if $u\in\mathcal{D}^{1,r}_0(\R^N) $ then $u$ vanishes at infinity, meaning that the set $\{x\in\R^N:|u(x)|\geq\epsilon\}$ has finite measure for any $\epsilon > 0$. In fact, by Chebichev's inequality and \eqref{embedding}, one has
\begin{equation*}
|\{x \in \R^N: \, |u(x)| \geq\epsilon\}| \leq \epsilon^{-r^*} \|u\|_{r^*}^{r^*} \leq (c \epsilon^{-1} \|u\|_{1,r})^{r^*} < +\infty,
\end{equation*}
where $ c > 0 $ is the best constant related to \eqref{embedding}. The monographs \cite{Ga,LL,SS} provide an exhaustive introduction on this topic. 

Finally,
\begin{equation*}
a\vee b:=\max\{a,b\},\;\; a\wedge b:=\min\{a,b\}\quad\forall\, a,b\in\R.
\end{equation*}
If $\Omega$ is a bounded domain in $\R^N$ then
\begin{equation*}
d(x):={\rm dist}(x,\partial\Omega),\quad x\in\overline{\Omega},   
\end{equation*}
while $\Omega'\subset\subset\Omega$ means $\overline{\Omega'}\subseteq \Omega$.
\section{Problems in bounded domains}\label{bound}
\subsection{The case $p=2$}\label{bound.1}
As far as we know, the study of singular semi-linear elliptic systems in bounded domains started with the paper \cite{CMcK1}, devoted to Gierer-Meinhardt's type \cite{GiMe} problem
\begin{equation}\label{modelone}
\left\{
\begin{alignedat}{2}
-\Delta u =-u+\frac{u}{v}\;\;&\mbox{in $\Omega$,}\;\; u>0\;\; && \mbox{in $\Omega$,}\;\; u\lfloor_{\partial\Omega}=0,\\
-\Delta v =-\alpha v+\frac{u}{v}\;\; & \mbox{in $\Omega$,}\;\; v>0\;\; && \mbox{in $\Omega$,}\;\; v\lfloor_{\partial\Omega}=0,
\end{alignedat}
\right.
\end{equation}
where $\Omega$ denotes a smooth bounded domain in $\R^N$, $N\geq 1$, while $\alpha\in\R^+$. More general right-hand sides were then considered in \cite{CMcK2,Kim,GR}. The work \cite{MS} investigates the system
\begin{equation}\label{modeltwo}
\left\{
\begin{alignedat}{2}
-\Delta u =\lambda u^{q_1}-\frac{u^{p_1}}{v^{\beta_1}}\;\;&\mbox{in $\Omega$,}\;\; u>0\;\; && \mbox{in $\Omega$,}\;\; u\lfloor_{\partial\Omega}=0,\\
-\Delta v =\mu v^{q_2}-\frac{u^{p_2}}{v^{\beta_2}}\;\; & \mbox{in $\Omega$,}\;\; v>0\;\; && \mbox{in $\Omega$,}\;\; v\lfloor_{\partial\Omega}=0.
\end{alignedat}
\right.
\end{equation}
Here, $q_1,q_2,\beta_1,\beta_2\in ]0,1[$, $p_1,p_2\in\R^+$, and $\lambda,\mu$ are two positive parameters. Define
\begin{equation*}
\sigma:=\frac{p_2(1-q_2)}{(1+\beta_2)(1-q_1)},\quad\eta:=\frac{\beta_1(1-q_1)}{(1-p_1)(1-q_2)}.   
\end{equation*}
An adequate sub-super-solution methods yields the following
\begin{thm}[\cite{MS}, Theorem 1.1]\label{thm1sem}
If $q_1<p_1$ then there exists $c_1>0$ such that \eqref{modeltwo} admits a solution $(u,v)\in C^{1,\alpha}(\overline{\Omega})^2$ for all $\lambda>0$, $\mu\geq c_1 \lambda^\sigma$.\\
If $q_1\geq p_1$ then there exists $c_2>0$ such that \eqref{modeltwo} has no solution once $\mu>0$ and $\lambda<c_2\mu^{-\eta}$.
\end{thm}
Two years before, in 2008, Hernandez, Mancebo, and Vega \cite{HMV} established the existence of classical solutions to the problem
\begin{equation}\label{modelthree}
\left\{
\begin{alignedat}{2}
\mathcal{L}_1 u =f(x,u,v)\;\;&\mbox{in $\Omega$,}\;\; u>0\;\;&& \mbox{in $\Omega$,}\;\; u\lfloor_{\partial\Omega}=0,\\
\mathcal{L}_2 v =g(x,u,v)\;\;&\mbox{in $\Omega$,}\;\; v>0\;\;&& \mbox{in $\Omega$,}\;\; v\lfloor_{\partial\Omega}=0,
\end{alignedat}
\right.
\end{equation}
where $\mathcal{L}_i$ denotes a linear, second-order, uniformly elliptic operator in non-divergence form while $f,g:\Omega\times(\R^+)^2\to\R$ are smooth enough. They examined both the cooperative and the non-cooperative situations; see \cite[Theorems 2.1--2.2]{HMV}. Here, cooperative means that the reaction $f(x,\cdot,\cdot)$ and $g(x,\cdot,\cdot)$ turn out increasing on $\R^+$ with respect to each variable separately. A \textit{uniqueness result} is also obtained provided \eqref{modelthree} turns out cooperative and concave, namely
\begin{equation*}
f(x,\tau s,\tau t)>\tau f(x,s,t),\;\;
g(x,\tau s,\tau t)>\tau g(x,s,t)
\end{equation*}
for all $\tau\in]0,1[$ and $(x,s,t)\in\Omega\times(\R^+)^2$; cf. \cite[Theorem 2.3]{HMV}. Many special cases are finally discussed. Singular systems driven by other linear, second-order, elliptic operators in divergence form were treated in \cite{BO}.

Let us next point out that the paper \cite{HS} contains a uniqueness result (cf. \cite[Theorem 2.2]{HS}) where no cooperative structure is supposed.

In 2015, Ghergu \cite{Ghe} thoroughly investigated both existence and uniqueness of classical solutions to the model system
\begin{equation}\label{modelfour}
\left\{
\begin{alignedat}{2}
-\Delta u=u^{\alpha_1}+v^{\beta_1}\;\;\mbox{in $\Omega$,}\;\; u>0\;\;\mbox{in $\Omega$,}\;\; && u\lfloor_{\partial\Omega}=0,\\
-\Delta v=v^{\alpha_2}+u^{\beta_2}\;\;\mbox{in $\Omega$,}\;\; v>0\;\;\mbox{in $\Omega$,}\;\; && v\lfloor_{\partial\Omega}=0
\end{alignedat}
\right.
\end{equation}
for $\alpha_i\vee\beta_1<0$, $i=1,2$. The work \cite{Mi} (see also \cite{Ghe1}) contains a similar analysis on the problem
\begin{equation}\label{modelfive}
\left\{
\begin{alignedat}{2}
-\Delta u=K_1(x) u^{\alpha_1}v^{\beta_1}\;\;\mbox{in $\Omega$,}\;\; u>0\;\;\mbox{in $\Omega$,}\;\; && u\lfloor_{\partial\Omega}=0,\\
-\Delta v=K_2(x)v^{\alpha_2}u^{\beta_2}\;\;\mbox{in $\Omega$,}\;\; v>0\;\;\mbox{in $\Omega$,}\;\; && v\lfloor_{\partial\Omega}=0,
\end{alignedat}
\right.
\end{equation}
where, as before, $\alpha_i\vee\beta_i<0$ and $K_i\in C^\alpha(\Omega)$, $i=1,2$.

The existence of solutions to singular convective elliptic systems was firstly studied by Alves and Moussaoui \cite{AM} in 2014.
\begin{thm}[\cite{AM}, Sections 3--4]
Let $-\alpha_i,\beta_i\in [0,1[$ and $g_i\in C^0(\R^{2N},\R^+)$ be bounded, $i=1,2$. Then the problem
\begin{equation*}
\left\{
\begin{alignedat}{2}
-\Delta u =v^{\alpha_1}\pm v^{\beta_1}+g_1(\nabla u,\nabla v)\;\;
& \mbox{in $\Omega$,}\;\; u>0\;\; && \mbox{in $\Omega$,}\;\; u\lfloor_{\partial\Omega}=0,\\
-\Delta v =u^{\alpha_2}\pm u^{\beta_2}+g_2(\nabla u,\nabla v)\;\; 
& \mbox{in $\Omega$,}\;\; v>0\;\; && \mbox{in $\Omega$,}\;\; v\lfloor_{\partial\Omega}=0
\end{alignedat}
\right.
\end{equation*}
admits a solution $(u,v)\in H^1_0(\Omega)^2\cap C^2(\Omega)^2$.
\end{thm}
Finally, the very recent work \cite{Calt} treats quasi-linear Schrödinger elliptic problems wiht both singular and convective reactions, while \cite{CMS} concerns singular systems having quadratic gradient. 
\subsection{Existence and multiplicity}\label{S2.1}
To shorten notation, given $p,q\in]1,+\infty[$, write
\begin{equation*}
\begin{split}
X^{p,q}(\Omega) &:=W^{1,p}(\Omega)\times W^{1,q}(\Omega),\quad
X^{p,q}_{\rm loc}(\Omega):=W^{1,p}_{\rm loc}(\Omega)\times W^{1,q}_{\rm loc}(\Omega),\\
& X^{p,q}_0(\Omega):=W^{1,p}_0(\Omega)\times W^{1,q}_0(\Omega).
\end{split}
\end{equation*}
From what we know, the first existence result for singular quasi-linear problems dates back to 2007 and examines the case
\begin{equation}\label{prob0}
\left\{
\begin{alignedat}{2}
-\Delta_p u & =v^{\alpha_1}+v^{\beta_1}\;\;\mbox{in $\Omega$,}\;\;
&& u>0\;\;\mbox{in $\Omega$,}\;\; u\lfloor_{\partial\Omega}=0,\\
-\Delta_q v & =u^{\alpha_2}+u^{\beta_2}\;\;\mbox{in $\Omega$,}\;\;
&& v>0\;\;\mbox{in $\Omega$,}\;\; v\lfloor_{\partial\Omega}=0.
\end{alignedat}
\right.
\end{equation}
Here, $\alpha_i<0<\beta_i$, while $N\geq 2$.
\begin{thm}[\cite{AC}, Theorem 1.1]\label{thm0}
If $\frac{2N}{N+1}\leq p,q<N$ and $0<-\alpha_i,\beta_i<\theta_i$, $i=1,2$, with
\begin{equation*}
\theta_1:=\min\left\{1,p-1,\frac{q'}{p'}\right\},\quad
\theta_2:=\min\left\{1,q-1,\frac{p'}{q'}\right\},
\end{equation*}
then \eqref{prob0} has a weak solution $(u,v)\in X^{p,q}_0(\Omega)$.  
\end{thm}
Its proof employs a result on nonlinear eigenvalue problems with lack of bifurcation due to Rabinowitz \cite{Ra} and a Hardy–Sobolev type inequality \cite{Ka}.

Now, consider the general problem
\begin{equation}\label{prob1}
\left\{
\begin{alignedat}{2}
-\Delta_p u & =f(x,u,v)\;\;\mbox{in $\Omega$,}\;\; && u>0\;\;\mbox{in $\Omega$,}\;\; u\lfloor_{\partial\Omega}=0,\\
-\Delta_q v & =g(x,u,v)\;\;\mbox{in $\Omega$,}\,\; && v>0\;\;\mbox{in $\Omega$,}\;\; v\lfloor_{\partial\Omega}=0,\\
\end{alignedat}
\right.
\end{equation}
where $f,g:\Omega\times(\R^+)^2\to\R$ fulfill Carath\'eodory's conditions. 
In 2012, El Manouni, Perera, and Shivaji investigated \eqref{prob1} under the assumptions below.
\begin{itemize}
\item[$({\rm a}_1)$] The functions $t\mapsto f(x,s,t)$ and $s\mapsto g(x,s,t)$ are increasing in $\R^+$, namely system \eqref{prob1} is \textit{cooperative}. 
\item[$({\rm a}_2)$] For every $0<s_0\leq s_1$ and $0<t_0\leq t_1$ one has
$$\sup_{\Omega\times[s_0,s_1]\times]0,t_1]}f(x,s,t)<+\infty,\;\;
\sup_{\Omega\times]0,s_1]\times[t_0,t_1]}g(x,s,t)<+\infty,$$
as well as
$$\sup_{\Omega\times[s_0,s_1]\times[t_0,t_1]}\max\left\{|f(x,s,t)|,
|g(x,s,t)|\right\}<+\infty.$$
\end{itemize}
They seek solutions $(u,v)\in X^{p,q}_{\rm loc}(\Omega)$ that satisfy the differential equations in the sense of distributions, i.e.,
\begin{equation*}
\begin{split}
\int_\Omega |\nabla u|^{p-2}\nabla u\cdot\nabla\varphi\,{\rm d}x=\int_\Omega f(\cdot, u,v)\varphi\,{\rm d}x\quad\forall\,\varphi\in C^\infty_0(\Omega),\\
\int_\Omega |\nabla v|^{q-2}\nabla v\cdot\nabla\psi\,{\rm d}x=\int_\Omega g(\cdot, u,v)\psi\,{\rm d}x\quad\forall\,\psi\in C^\infty_0(\Omega),
\end{split}
\end{equation*}
and, moreover, $u,v\in C^0(\overline{\Omega})$.

A pair $(\underline{u},\underline{v})\in X^{p,q}(\Omega)$ is called a sub-solution of problem \eqref{prob1} provided $f(\cdot,\underline{u},\underline{v}) \in L^{p'}(\Omega)$, $g(\cdot,\underline{u},\underline{v})\in L^{q'}(\Omega)$, and
\begin{equation*}
-\Delta_p\underline{u}\leq f(\cdot,\underline{u}, \underline{v})\;\;\mbox{in $\Omega$,}\;\;
-\Delta_p\underline{v}\leq g(\cdot,\underline{u}, \underline{v})\;\;\mbox{in $\Omega$,}\;\;
\underline{u}\vee\underline{v}\leq 0\;\:\mbox{on $\partial\Omega$.}
\end{equation*}
A super-solution $(\overline{u},\overline{v})$ is defined similarly, by reversing all the above inequalities.

Let $\{\epsilon_n\}\subseteq\R^+$ satisfy $\epsilon_n\to 0$. Set, for every $n\in\N$, $(x,s,t)\in\Omega\times(\R^+)^2$,
\begin{equation*}
f_n(x,s,t):=f(x,s\vee\epsilon_n,t\vee\epsilon_n),\;\;
g_j(x,s,t):=g(x,s\vee\epsilon_n,t\vee\epsilon_n),
\end{equation*}
and consider the sequence of regularized systems
\begin{equation}\label{prob1aux}
\left\{
\begin{alignedat}{2}
-\Delta_p u & =f_n(x,u,v)\;\; &&\mbox{in $\Omega$,} \;\; u\lfloor_\Omega=0,\\
-\Delta_q v & =g_n(x,u,v)\;\; &&\mbox{in $\Omega$,} \;\; v\lfloor_\Omega=0.
\end{alignedat}
\right.
\end{equation}
\begin{thm}[\cite{EPS}, Theorem 3.1]\label{thm1}
Suppose $({\rm a}_1)$--$({\rm a}_2)$ hold. If, for each $n\in\N$, there exist a sub-solution $(\underline{u}_n,\underline{v}_n)$ and a super-solution $(\overline{u}_n,\overline{v}_n)$ to \eqref{prob1aux} such that $(\underline{u}_n,\underline{v}_n)\leq (\overline{u}_n,\overline{v}_n)$,
\begin{equation*}
\inf_{n\in\N}\essinf_{\Omega'}\,(\underline{u}_n\wedge\underline{v}_n)>0
\end{equation*}
whenever $\Omega'\subset\subset\Omega$, and
\begin{equation*}
\sup_{n\in\N}\esssup_{\Omega}\,(\overline{u}_n\vee\overline{v}_n)<+\infty,
\end{equation*}
then \eqref{prob1} admits a distributional solution $(u,v)\in C^{1,\alpha}_{\rm loc}(\Omega)^2 \cap C^0(\overline{\Omega})^2$.
\end{thm}
Sufficient conditions for the existence of sub-super-solution pairs to \eqref{prob1aux} are given in \cite[Sections 4--5]{EPS}. As an example, via Theorem \ref{thm1} one can show that the model problem (cf. \eqref{modelfour})
\begin{equation*}
\left\{
\begin{alignedat}{2}
-\Delta_p u=u^{\alpha_1}+\mu v^{\beta_1}\;\;\mbox{in $\Omega$,}\;\; u>0\;\;\mbox{in $\Omega$,}\;\; && u\lfloor_{\partial\Omega}=0,\\
-\Delta_q v=v^{\alpha_2}+\mu u^{\beta_2}\;\;\mbox{in $\Omega$,}\;\; v>0\;\;\mbox{in $\Omega$,}\;\; && v\lfloor_{\partial\Omega}=0,
\end{alignedat}
\right.
\end{equation*}
where $\alpha_1\vee\alpha_2<0$ and $\beta_1\wedge\beta_2\geq 0$, possesses a solution provided $\mu\geq 0$ is small enough. Moreover, singular \textit{semi-positone systems}, which means both $\lim_{s\to 0^+}f(x,s,t)=-\infty$ uniformly in $(x,t)$ and $\lim_{t\to 0^+}g(x,s,t)=-\infty$ uniformly with respect to $(x,s)$, are investigated. Let us also mention the papers \cite{LSY}, \cite{KM}, and \cite{CHS}. In particular, \cite{CHS} deals with a non-cooperative system.

One year later, Giacomoni, Hernandez, and Sauvy obtained the general results below, where $f,g\in C^1(\Omega\times(\R^+)^2)$, through a different notion of sub-super-solution.

We say that $(\underline{u},\underline{v}),(\overline{u},\overline{v})\in X^{p.q}_{\rm loc}(\Omega)\cap C^0(\overline{\Omega})^2$ are a sub-super-solution pair to \eqref{prob1} when $\underline{u}$, $\underline{v}$, $\overline{u}$, $\overline{v}$ are 
locally uniformly positive,
$(\underline{u},\underline{v})\leq (\overline{u},\overline{v})$, and for any $(u,v)\in [\underline{u},\overline{u}]\times [\underline{v},\overline{v}]$ one has 
\begin{equation*}
-\Delta_p\underline{u}\leq f(\cdot,\underline{u},v),\;
-\Delta_p\underline{v}\leq g(\cdot,u,\underline{v}),\;
-\Delta_p\overline{u}\geq f(\cdot,\overline{u}, v),\;
-\Delta_p\overline{v}\geq g(\cdot,u,\overline{v})
\end{equation*}
in $\Omega$. The set ${\cal C}:=[\underline{v},\overline{v}]\times [\underline{u},\overline{u}]$ is usually called \textit{trapping region}; cf. \cite{CaMo}.
\begin{thm}[\cite{GHS}, Theorem 2.1]\label{thm2}
Let $(\underline{u},\underline{v}),(\overline{u},\overline{v})\in X^{p,q}_0(\Omega)$ be a sub-super-solution pair to \eqref{prob1} fulfilling:
\begin{itemize}
\item[$({\rm a}_3)$] $\overline {u}\leq c_1d^{\gamma_1}$ and $\overline{v}\leq c_2d^{\gamma_2}$, with appropriate $c_1,c_2,\gamma_1,\gamma_2>0$.
\item[$({\rm a}_4)$] There exist $c_3,c_4>0$, $\delta_1,\delta_2 \in\R$ such that 
\begin{equation*}
|f(\cdot,u,v)|\leq c_3 d^{\delta_1},\,\;|g(\cdot,u,v)|\leq c_4 d^{\delta_2}\,\;\forall\, (u,v)\in{\cal C}.
\end{equation*}
\item[$({\rm a}_5)$] For suitable $c_5,c_6,\sigma_1,\sigma_2>0$ one has
\begin{equation*}
\left|\frac{\partial f}{\partial s}(\cdot,u,v)\right|\leq c_5d^{\delta_1-\sigma_1},\;
\left|\frac{\partial g}{\partial t}(\cdot,u,v)\right|\leq c_6d^{\delta_2-\sigma_2},
\;(u,v)\in {\cal C}.
\end{equation*}
\end{itemize}
Assume further that
\begin{equation*}
\delta_1>-2+\frac{1}{p}+(\sigma_1-\gamma_1)^+,\;\;
\delta_2>-2+\frac{1}{q}+(\sigma_2-\gamma_2)^+.
\end{equation*}
Then \eqref{prob1} admits a weak solution $(u,v)\in{\cal C}$.
\end{thm}
\begin{thm}[\cite{GHS}, Theorem 2.3]\label{thm3}
Suppose $(\underline{u},\underline{v}),(\overline{u},\overline{v})$ is a sub-super-solution pair to \eqref{prob1} complying with $({\rm a}_3)$. If
\begin{equation}\label{hypghs1}
\frac{\partial f}{\partial t}(x,s,t)>0,\;\;\frac{\partial g}{\partial s}(x,s,t)>0\quad\forall\, (x,s,t)\in\Omega\times(\R^+)^2
\end{equation}
and there exist $c_7,c_8>0$, $\eta_1,\eta_2\in\R$ such that
\begin{equation*}
\left|\frac{\partial f}{\partial s}(\cdot,u,v)\right|\leq c_7d^{\eta_1},\;\;
\left|\frac{\partial g}{\partial t}(\cdot,u,v)\right|\leq c_8d^{\eta_2}
\;\;\forall\, (u,v)\in {\cal C},
\end{equation*}
then \eqref{prob1} possesses a distributional solution $(u,v)\in{\cal C}$.
\end{thm}
Proofs are based on a very nice, non-trivial use of Schauder's fixed point theorem. Applications to the model problem 
\begin{equation}\label{prob1aux2}
\left\{
\begin{alignedat}{2}
-\Delta_p u =K_1(x) u^{\alpha_1}\, v^{\beta_1}\;\; & \mbox{in $\Omega$,}\;\; u>0\;\; && \mbox{in $\Omega$,}\;\; u\lfloor_{\partial\Omega}=0,\\
-\Delta_q v =K_2(x) v^{\alpha_2}\, u^{\beta_2}\;\; & \mbox{in $\Omega$,}\;\;
v>0\;\; && \mbox{in $\Omega$,}\;\; v\lfloor_{\partial\Omega}=0,
\end{alignedat}
\right.
\end{equation}
where $K_1,K_2$ satisfy appropriate conditions while
\begin{equation}\label{hypghs}
\alpha_1<p-1,\;\;\alpha_2<q-1,\;\;
(p-1-\alpha_1)(q-1-\alpha_2)>|\beta_1\beta_2|>0,
\end{equation}
are given. Evidently, \eqref{prob1aux2} becomes \eqref{modelfive} for $p=q=2$. See also \cite{DOS}, where $p=q<N$, $K_1\in L^\infty(\Omega)_+$, $K_2\in L^\gamma(\Omega)_+$ for some $\gamma>\frac{N}{p}$, 
\begin{equation*}
0<\beta_1<1\wedge (p-1),\;\;-1<\alpha_1<p-\beta_1-1,\;\;\alpha_2:=\beta_1-1,\;\;
\beta_2:=\alpha_1+1\, .
\end{equation*}

A special case of \eqref{prob1aux2} was treated in \cite{MM1} (cf. in addition \cite{DeC}), namely
\begin{equation}\label{prob1aux3}
\left\{
\begin{alignedat}{2}
-\Delta_p u =u^{\alpha_1}\, v^{\beta_1}\;\; & \mbox{in $\Omega$,}\;\; u>0\;\; && \mbox{in $\Omega$,}\;\; u\lfloor_{\partial\Omega}=0,\\
-\Delta_q v =v^{\alpha_2}\, u^{\beta_2}\;\; & \mbox{in $\Omega$,}\;\;
v>0\;\; && \mbox{in $\Omega$,}\;\; v\lfloor_{\partial\Omega}=0.
\end{alignedat}
\right.
\end{equation}
\begin{thm}[\cite{MM1}, Theorem 1.1]\label{thm5}
Let $1<p,q\leq N$ and let $\alpha_i,\beta_i$ satisfy
\begin{equation}\label{hypmm1}
\begin{split}
-2+\frac{1}{p}\leq\alpha_1<0,\;\;-2+\frac{1}{q}\leq\alpha_2<0,\\
0<\beta_1<\frac{q}{p}(p-1-\alpha_1),\;\;0<\beta_2<\frac{p}{q}(q-1-\alpha_2).
\end{split}
\end{equation}
Then \eqref{prob1aux3} has a weak solution $(u,v)\in X^{p,q}_0(\Omega) \cap L^\infty(\Omega)^2$.
\end{thm}
\begin{rmk}
It should be noted that \eqref{hypmm1} forces \eqref{hypghs}. Moreover, both \eqref{hypghs1} and \eqref{hypmm1} basically entail $({\rm a}_1)$. 
\end{rmk}
System \eqref{prob1aux3}, with a \textit{competitive} interaction between the two components $u$ and $v$ was thoroughly studied by Giacomoni, Schindler, and Takac \cite{GST}. Recall that \eqref{prob1} is said to be competitive if
\begin{itemize}
\item[$({\rm a}'_1)$] The functions $t\mapsto f(x,s,t)$ and $s\mapsto g(x,s,t)$ are decreasing on $\R^+$.
\end{itemize}
In the situation above this means $\alpha_i\vee\beta_i<0$. Under the \textit{sub-homogeneity condition}
\begin{equation}\label{subhom}
 (p-1-\alpha_1)(q-1-\alpha_2)>\beta_1\beta_2   
\end{equation}
and suitable upper bounds regarding $\alpha_i,\beta_i$, they proved that \eqref{prob1aux3} admits a solution $(u,v)\in X^{p,q}_0(\Omega)\cap C^{0,\alpha}(\overline{\Omega})^2$. When $p=q$, see also \cite{Si}.

Existence results for problem \eqref{prob1} where the competitive structure $({\rm a}'_1)$ is allowed can be found in \cite{MM2,MM3}.

The recent work \cite{AMT} examines singular $(p(x),q(x))$-Laplacian problems with singularity coming through logarithmic reactions that involve variable exponents growth conditions; cf. also \cite{AM1,MV}.

As far as we know, till today, much less attention has been paid to \textit{multiplicity of solutions}. Actually, we can only mention the papers \cite{S,CSGS,DM,AN,KS}. The first deals with singular $p(x)$-Laplacian systems while the second is devoted to quasi-linear problems driven by $(\Phi_1,\Phi_2)$-Laplace operators. Theorem 1 in \cite{DM} considers the case when $f,g$ do not depend on $x$, are positive, entail a cooperative structure, and, roughly speaking, cross the first eigenvalue at infinity, namely
\begin{equation*}
\lim_{s\to+\infty}\frac{f(s,t)}{s^{p-1}}=J_1>\lambda_{1,p},\quad
\lim_{t\to+\infty}\frac{g(s,t)}{t^{q-1}}=J_2>\lambda_{1,q}.
\end{equation*}
Two smooth solutions are obtained combining  sub-super-solution methods with the Leray-Schauder topological degree.
A different approach is adopted in \cite{AN}. The differential operators, which include the $r$-Laplacian as a special case, are neither homogeneous nor linear, while
\begin{equation*}
f(x,s,t):=K_1(x)s^{\alpha_1}+\frac{\partial h}{\partial s}(x,s,t),\quad
g(x,s,t):=K_2(x)t^{\alpha_2}+\frac{\partial h}{\partial t}(x,s,t),
\end{equation*}
where $K_i\in L^\infty(\Omega)$, $\alpha_i<0$, and $h\in C^1(\overline{\Omega} \times\R^2)$ satisfy appropriate conditions that permit the use of variational methods. Theorem 1.7 of \cite{AN} gives two positive solutions provided $\max_{i=1,2}\Vert K_i\Vert_\infty$ is sufficiently small. Finally, \cite{KS} addresses singular $p(x)$-Laplacian systems with nonlinear boundary conditions.
\subsection{Uniqueness}
Unless the semi-linear case, uniqueness of solutions looks a difficult matter, even for problem \eqref{prob1aux3}. In fact, as pointed out by Giacomoni, Hernandez, and Moussaoui \cite{GHM}, equations with quasi-linear elliptic operators exhibit additional troubles for obtaining the validity of the strong comparison principle, which requires the $C^1$-regularity of solutions. If it cannot be obtained (as in the strongly singular setting $\alpha_i+\beta_i<-1$) then one can still try to use a suitable variant of the well known Krasnoleskii’s argument \cite{K}.\\
Let us first examine the cooperative case $\beta_1\wedge\beta_2>0$. Proposition 3.1 of \cite{GHM} basically yields
\begin{thm}[\cite{GHM}, Theorem 3.2]\label{thm6}
If $-1<\alpha_i+\beta_i<0<\beta_i$, $i=1,2$, then \eqref{prob1aux3} admits a unique weak solution $(u,v)\in X^{p,q}_0(\Omega)\cap C^{1,\alpha} (\overline{\Omega})^2$. Moreover, $u,v\geq c_1d$ for some $c_1>0$.
\end{thm}
Suppose now
\begin{equation}\label{strongsing1}
 \alpha_1+\beta_i<-1,\; i=1,2,\;\;0<\beta_1<p-1,\;\;0<\beta_2<q-1,
\end{equation}
and denote by $\gamma,\theta\in]0,1[$ the unique solution to the system
\begin{equation}\label{strongsing2}
\left\{
\begin{aligned}
(\gamma-1)(p-1)-1 & =\alpha_1\gamma+\beta_1\theta\\   
(\theta-1)(q-1)-1 & =\beta_2\gamma+\alpha_2\theta.
\end{aligned}
\right.
\end{equation}
\begin{thm}[\cite{GHM}, Theorem 3.3]\label{thm7}
Let \eqref{strongsing1}--\eqref{strongsing2} be satisfied. Assume also that
\begin{equation}\label{strongsing3}
\alpha_1\gamma+\beta_1\theta>\frac{1}{p}-2,\;\; \beta_2\gamma+\alpha_2\theta>\frac{1}{q}-2 .  
\end{equation}
Then \eqref{prob1aux3} possesses a unique weak solution $(u,v)\in X^{p,q}_0(\Omega)$ fulfilling
\begin{equation*}
c_2(\varphi_{1,p}^\gamma,\varphi_{1,q}^\theta)\leq (u,v)\leq
c_3(\varphi_{1,p}^\gamma,\varphi_{1,q}^\theta)
\end{equation*}
with appropriate $c_2,c_3>0$.
\end{thm}
We next examine the competitive case $\beta_i<0$. Via the sub-homogeneity condition \eqref{subhom}, Theorem 2.2 of \cite{GST} considers various possible choices of exponents. Here, for the sake of brevity, we will present only one.
\begin{thm}[\cite{GST}, Theorem 2.2]\label{thm8}
Suppose $\alpha_i\vee\beta_i<0$. Under \eqref{subhom}, \eqref{strongsing2}, and \eqref{strongsing3}, the same conclusion of Theorem \ref{thm7} is true.
\end{thm}
\subsection{Systems with convection terms}
In 2017, Motreanu, Moussaoui, and Zhang \cite{MMZ} treated the general problem
\begin{equation}\label{probConv}
\left\{
\begin{alignedat}{2}
-\Delta_p u =f(u,v,\nabla u,\nabla v)\;\; & \mbox{in $\Omega$,}\;\; u>0\;\; && \mbox{in $\Omega$,}\;\; u\lfloor_{\partial\Omega}=0,\\
-\Delta_q v =g(u,v,\nabla u,\nabla v)\;\; & \mbox{in $\Omega$,}\;\; v>0\;\; && \mbox{in $\Omega$,}\;\; v\lfloor_{\partial\Omega}=0,
\end{alignedat}
\right.
\end{equation}
where $f,g:(\R^+)^2\times(\R^N\setminus\{0\})^2\to\R^+$ are continuous functions satisfying the growth conditions below.
\begin{itemize}
\item[$({\rm a}_6)$] There exist $\mu_i,\hat\mu_i>0$ and $\alpha_i,\beta_i\in\R$, such that
\begin{equation*}
\begin{aligned}
\mu_1 (s+|\xi_1|)^{\alpha_1}(t+|\xi_2|)^{\beta_1}\leq f(s,t,\xi_1, \xi_2)
\leq\hat\mu_1 (s+|\xi_1|)^{\alpha_1}(t+|\xi_2|)^{\beta_1}\\
\mu_2 (s+|\xi_1|)^{\beta_2}(t+|\xi_2|)^{\alpha_2}\leq g(s,t,\xi_1,\xi_2)\leq \hat\mu_2 (s+|\xi_1|)^{\beta_2}(t+|\xi_2|)^{\alpha_2}
\end{aligned}
\end{equation*}
for all $(s,t,\xi_1,\xi_2)\in (\R^+)^2\times(\R^N\setminus\{0\})^2$.
\item[$({\rm a}_7)$] One has
\begin{equation*}
\beta_1\vee\beta_2<0,\quad\alpha_1\alpha_2>0,\quad |\alpha_1|<p-1+\beta_1,\quad |\alpha_2|<q-1+\beta_2.
\end{equation*}
\end{itemize}
Since $\beta_1\vee\beta_2<0$, system \eqref{probConv} has a competitive structure. Moreover, \textit{the right-hand sides may exhibit singularities in both the solution and its gradient.}

Combining comparison arguments with a priori estimates yields the next
\begin{thm}[\cite{MMZ}, Theorem 1.1]\label{thm9}
If $({\rm a}_6)$--$({\rm a}_7)$ hold then \eqref{probConv} admits a solution $(u,v)\in C^{1,\alpha}(\overline{\Omega})^2$.
\end{thm}
Three years later, also the cooperative case was examined; see \cite{CLM}. Now, $1<p,q<N$ while $f,g:(\R^+)^2\times\R^{2N}\to\R^+$ are continuous and fulfill 
\begin{itemize}
\item[$({\rm a}_8)$] There exist $\mu_i,\hat\mu_i,\beta_i,\gamma_i,\theta_i\in\R^+$, $\alpha_i\in]-1,0[$, $i=1,2$, such that
\begin{equation*}
\begin{split}
\mu_1 s^{\alpha_1} t^{\beta_1}\leq f(s,t,\xi_1, \xi_2)\leq
\hat\mu_1 s^{\alpha_1} t^{\beta_1}+|\xi_1|^{\gamma_1}+|\xi_2|^{\theta_1}, \\
\mu_2 s^{\beta_2} t^{\alpha_2}\leq g(s,t,\xi_1, \xi_2)\leq
\hat\mu_2 s^{\beta_2} t^{\alpha_2}+|\xi_1|^{\gamma_2}+|\xi_2|^{\theta_2}
\end{split}
\end{equation*}
for all $(s,t,\xi_1,\xi_2)\in(\R^+)^2\times\R^{2N}$.
\item[$({\rm a}_9)$] One has $\alpha_i+\beta_i\geq 0$, as well as
\begin{equation*}
\max\{-\alpha_1+\beta_1,\gamma_1,\theta_1\}<p-1,\;\;
\max\{-\alpha_2+\beta_2,\gamma_2,\theta_2\}<q-1.
\end{equation*}
\end{itemize}
\begin{thm}[\cite{CLM}, Theorem 1]\label{thm10}
Under $({\rm a}_8)$--$({\rm a}_9)$, problem \eqref{probConv} possesses a solution $(u,v)\in C^1_0(\overline{\Omega})^2$. Moreover, $c_1d\leq u,v\leq c_2d$ for suitable $c_1,c_2>0$.
\end{thm}
A further interesting contribution in contained in \cite{DeMou}. The very recent paper \cite{GM} (see also \cite{GM1}) establishes the existence of infinitely many solutions to the Neumann problem
\begin{equation}\label{probNeumann}
\left\{
\begin{alignedat}{2}
-\Delta_p u =f(x,u,v,\nabla u,\nabla v)\;\; & \mbox{in $\Omega$,}\;\; u>0\;\; && \mbox{in $\Omega$,}\;\;
\frac{\partial u}{\partial\nu}=0\;\;\mbox{on $\partial\Omega$,}\\
-\Delta_q v =g(x,u,v,\nabla u,\nabla v)\;\; & \mbox{in $\Omega$,}\;\; v>0\;\; && \mbox{in $\Omega$,}\;\;
\frac{\partial v}{\partial\nu}=0\;\;\mbox{on $\partial\Omega$,}
\end{alignedat}
\right.
\end{equation}
where $1<p,q<+\infty$, $f,g:\Omega\times (\R^+)^2\times\R^{2N} \to\R$ satisfy Caratheodory's conditions while $\nu$ denotes the outer unit normal to $\partial\Omega$.\\
The \textit{sub-linear case} is first investigated using the next assumptions.
\begin{itemize}
\item[$({\rm a}_{10})$] There exist $\alpha_i<0<\beta_i$, $\gamma_1,\delta_1\in[0,p-1[$, $\gamma_2,\delta_2\in[0,q-1[$, and $a_i,b_i,c_i\in L^\infty(\Omega)$ such that
\begin{equation*}
\begin{aligned}
|f(x,s,t,\xi_1,\xi_2)| \leq a_1(x) s^{\alpha_1}t^{\beta_1} + b_1(x) (|\xi_1|^{\gamma_1} + |\xi_2|^{\delta_1}) + c_1(x),\\
|g(x,s,t,\xi_1,\xi_2)| \leq a_2(x) s^{\beta_2}t^{\alpha_2} + b_2(x) (|\xi_1|^{\gamma_2} + |\xi_2|^{\delta_2}) + c_2(x)
\end{aligned}
\end{equation*}
for all $(x,s,t,\xi_1,\xi_2)\in\Omega\times(\R^+)^2\times \R^{2N}$.

\item[$({\rm a}_{11})$] There are $\{h_n\},\{\hat{h}_n\}, \{k_n\},\{\hat{k}_n\}, \{C_n\}\subseteq\R^+$, with $ C_n \to +\infty $, satisfying $h_n<k_n<h_{n+1}$, $\hat{h}_n<\hat{k}_n<\hat{h}_{n+1}$,
\begin{equation}\label{None}
\begin{split}
f(x,k_n,t,\xi_1,\xi_2) &\leq 0 \leq f(x,h_n,t,\xi_1,\xi_2),\\
g(x,s,\hat{k}_n,\xi_1,\xi_2) &\leq 0 \leq g(x,s,\hat{h}_n,\xi_1,\xi_2)
\end{split}
\end{equation}
for all $ (x,s,t,\xi_1,\xi_2) \in \Omega \times [h_n,k_n] \times [\hat{h}_n,\hat{k}_n] \times B_{\R^N}(C_n)^2 $, $n\in\N$, and
\begin{equation}\label{Ntwo}
\lim_{n\to\infty}\frac{h_n^{\alpha_1}\hat{k}_n^{\beta_1}}{C_n^{p-1}}=0,\quad
\lim_{n\to\infty}\frac{\hat{h}_n^{\alpha_2}k_n^{\beta_2}}{C_n^{q-1}}=0.
\end{equation}
\end{itemize}
\begin{thm}[\cite{GM}, Theorem 4.2]
If $({\rm a}_{10})$--$({\rm a}_{11})$ hold then \eqref{probNeumann} has a sequence of solutions $\{(u_n,v_n)\}\subseteq C^1(\overline{\Omega})^2$ such that $(u_n,v_n)<(u_{n+1},v_{n+1})$ for every $n\in\N$. Moreover,
$\displaystyle{\lim_{n\to\infty}}u_n=\displaystyle{\lim_{n\to\infty}}v_n=+\infty$ uniformly in $\overline{\Omega}$ once $ h_n,\hat{h}_n \to +\infty $.
\end{thm}
As regards the \textit{super-linear case}, denote by $({\rm a}'_{10})$ condition $({\rm a}_{10})$ written for $\gamma_1,\delta_1>p-1$ and $\gamma_2,\delta_2>q-1$. Similarly,
\begin{itemize}
\item[$({\rm a}'_{11})$] There exist $\{h_n\},\{\hat{h}_n\},\{k_n\}, \{\hat{k}_n\},\{C_n\}\subseteq\R^+$, with $ C_n \to 0 $, satisfying
$k_{n+1}<h_n<k_n$, $\hat{k}_{n+1}<\hat{h}_n<\hat{k}_n$ for all $n\in\N$ and \eqref{None}--\eqref{Ntwo}.
\end{itemize}
\begin{thm}[\cite{GM}, Theorem 4.3]\label{superlinearcase}
Under $({\rm a'}_{10})$--$({\rm a'}_{11})$,  problem \eqref{probNeumann} possesses a sequence of solutions $\{(u_n,v_n)\}\subseteq C^1(\overline{\Omega})^2$ such that $(u_{n+1},v_{n+1})<(u_n,v_n)$ for every $n\in\N$. Moreover, $\displaystyle{\lim_{n\to\infty}}u_n=\displaystyle{\lim_{n\to\infty}}v_n=0$ uniformly in $\overline{\Omega}$ once $ k_n,\hat{k}_n \to 0 $.
\end{thm}
An easy example of nonlinearities, with both singular and convective terms, that fulfill \eqref{None}--\eqref{Ntwo} is the following.
\begin{ex}
Set, for every $(x,s,t,\xi_1,\xi_2)\in\Omega\times(\R^+)^2\times\R^{2N}$,
\begin{equation*}
\begin{split}
f(x,s,t,\xi_1,\xi_2) = \sin\frac{1}{s}\,\left(s^{\alpha_1} t^{\beta_1} - |\xi_1|^{\gamma_1} - |\xi_2|^{\delta_1}\right), \\
g(x,s,t,\xi_1,\xi_2) = \cos\frac{1}{t}\,\left(s^{\beta_2} t^{\alpha_2} - |\xi_1|^{\gamma_2} - |\xi_2|^{\delta_2}\right),
\end{split}
\end{equation*}
where
\begin{equation*}
\gamma_1\wedge\delta_1>\alpha_1+\beta_1 > p-1,\quad
\gamma_2\wedge\delta_2> \alpha_2+\beta_2 > q-1. 
\end{equation*}
To check \eqref{None}--\eqref{Ntwo} one can pick $C_n=\frac{1}{n}$,
\begin{equation*}
h_n=\frac{1}{\pi/2 +2\pi n},\; k_n=\frac{1}{-\pi/2+2\pi n},\; \hat{h}_n=\frac{1}{2\pi+2\pi n},\;\hat{k}_n=\frac{1}{\pi+2\pi n}.
\end{equation*}
\end{ex}
\section{Problems on the whole space}\label{wholespace}
\subsection{The case $p=2$}
In 2009, Moussaoui, Khodja, and Tas \cite{MKT} studied the following singular,, semi-linear elliptic, Gierer-Meinhardt's type system:
\begin{equation}\label{wholemodelone}
\left\{
\begin{alignedat}{2}
-\Delta u+\alpha_1(x) u=a_1(x)\frac{1}{v^q}\;\; & \mbox{in $\R^N$,}\;\; u>0\;\; && \mbox{in $\R^N$,}\\
-\Delta v+\alpha_2(x) v=a_2(x)\frac{u^r}{v^s}\;\; & \mbox{in $\R^N$,}\;\; v>0\;\; && \mbox{in $\R^N$,}\\
u(x)\to 0,\; v(x)\to 0\;\; & \mbox{as $|x|\to\infty$,}
\end{alignedat}
\right.
\end{equation}
where, roughly speaking, $\alpha_i,a_i\in L^\infty_{\rm loc}(\R^N)_+$, $a_i$ satisfies suitable integrability conditions, $q,r,s>0$, and $r\leq s+1<2$. A solution $(u,v)\in \mathcal{D}^{1,2}_0(\R^N)^2$ is obtained via Schauder's fixed point theorem. Previous papers on the same subject are \cite{DKC,DKW}, whilst, excepting \cite{B1,B2}, we were not able to find more recent contributions.

Let us next point out that, even for the semi-linear case, the question whether a singular elliptic problem on the whole space admits multiple positive solutions is an open question.

Finally, as regards uniqueness, we mention \cite[Theorem 1.1]{B1}, which deals with the system
\begin{equation}\label{wholemodeltwo}
\left\{
\begin{alignedat}{2}
-\Delta u+\alpha(x) u^2=a(x)\frac{v^{1-p}}{u^p}\;\; & \mbox{in $\R^N$,}\;\; u>0\;\; && \mbox{in $\R^N$,}\\
-\Delta v+\alpha(x) v^2=a(x)\frac{u^{1-p}}{v^p}\;\; & \mbox{in $\R^N$,}\;\; v>0\;\; && \mbox{in $\R^N$,}
\end{alignedat}
\right.
\end{equation}
where $N\geq 3$, $p\in]\frac{1}{2},1[$, $\alpha\in L^1(\R^N)\cap L^{\frac{q}{q-2}}(\R^N)$ for some $q\in]2,2^*[$, $\alpha>0$, and $0\leq a\leq\alpha$.

Existence and uniqueness of solutions to singular convective elliptic problems was firstly studied by Benrhouma \cite{B2} in 2017, who considered the system
\begin{equation*}
\left\{
\begin{alignedat}{2}
-\Delta u+\alpha_1\frac{|\nabla u|^2}{u}=\frac{p}{p+q}a(x)u^{p-1} v^q+b_1(x)\;\; & \mbox{in $\R^N$,}\;\; u>0\;\; && \mbox{in $\R^N$,}\\
-\Delta v+\alpha_2\frac{|\nabla v|^2}{v}=\frac{q}{p+q}a(x)u^p v^{q-1}+b_2(x)\;\; & \mbox{in $\R^N$,}\;\; v>0\;\; && \mbox{in $\R^N$,}
\end{alignedat}
\right.
\end{equation*}
with $p\wedge q>1$, $p+q<2^*-1$, $a$ fulfilling appropriate integrability conditions, $\alpha_i>\frac{N+2}{4}$, and $b_i\in L^{2^*}(\R^N)_+\cap L^\infty(\R^N)$, $i=1,2$.
\subsection{Existence and multiplicity}
Henceforth, given $p,q\in]1,+\infty[$, we will write
\begin{equation*}
X^{p,q}(\R^N):=\mathcal{D}^{1,p}_0(\R^N)\times\mathcal{D}^{1,q}_0(\R^{N}).
\end{equation*}
To the best of our knowledge, until 2019, singular elliptic systems in the whole space were investigated only for $p:=q:=2$, essentially exploiting the linearity of involved differential operators. The paper \cite{MMM} considers the problem
\begin{equation}\label{probwholeone}
\left\{
\begin{alignedat}{2}
-\Delta_p u =a_1(x)f(u,v)\;\; & \mbox{in $\R^N$,}\;\; u>0\;\; && \mbox{in $\R^N$,}\\
-\Delta_q v =a_2(x)g(u,v)\;\; & \mbox{in $\R^N$,}\;\; v>0\;\; && \mbox{in $\R^N$,}
\end{alignedat}
\right.
\end{equation}
where $N\geq 3$ while $1<p,q<N$. Nonlinearities $f,g:(\R^+)^2\to \R^+$ are continuous and fulfill the condition
\begin{itemize}
\item[$({\rm a}_{12})$] There exist $\mu_i,\hat\mu_i>0$, $i=1,2$, such that
\begin{equation*}
\mu_1 s^{\alpha_1}\leq f(s,t)\leq\hat\mu_1 s^{\alpha_1} (1+t^{\beta_1}),\;\;
\mu_2 t^{\alpha_2}\leq g(s,t)\leq\hat\mu_2 (1+s^{\beta_2}) t^{\alpha_2}
\end{equation*}
for all $s,t\in\R^+$, with $-1<\alpha_i<0<\beta_i$,
\begin{equation}\label{c2}
\alpha_1+\beta_2<p-1,\;\;\alpha_2+\beta_1<q-1,
\end{equation}
as well as
\begin{equation*}
\beta_1<\frac{q^*}{p^*}\min\{p-1, p^*-p\},\;\;
\beta_2<\frac{p^*}{q^*}\min\{q-1, q^*-q\}.
\end{equation*}
\end{itemize}
Coefficients $a_i:\R^N\to\R^+$ satisfy the assumption
\begin{itemize}
\item[$({\rm a}_{13})$] $a_i\in L^1(\R^N)\cap L^{\zeta_i}(\R^N)$, where
\begin{equation*}
\frac{1}{\zeta_1}\leq 1-\frac{p}{p^*}-\frac{\beta_1}{q^*}\, ,\;\;\frac{1}{\zeta_2}\leq 1- \frac{q}{q^*}- 
\frac{\beta_2}{p^*}\, .
\end{equation*}
\end{itemize}
A pair $(u,v)\in X^{p,q}(\R^N)$ is called a (weak) solution to \eqref{probwholeone} provided $u,v>0$ and
\begin{equation*} 
\begin{aligned}
\int_{\R^N}|\nabla u|^{p-2}\nabla u\nabla \varphi\,{\rm d}x & =\int_{\R^N}a_{1}f(u,v)\varphi\,{\rm d}x\;\;\forall\,\varphi\in \mathcal{D}^{1,p}(\R^N) ,\\
\int_{\R^N}|\nabla v|^{q-2}\nabla v\nabla \psi\, dx & =\int_{\R^N}a_{2}g(u,v)\psi\, dx\;\;\forall\,\psi\in \mathcal{D}^{1,q}(\R^N).
\end{aligned}
\end{equation*}
Variational methods do not work, at least in a direct way, because the Euler functional associated with \eqref{probwholeone} is not well defined. A similar comment holds for sub-super-solution techniques, that are usually employed in the case of bounded domains. So, one is naturally led to apply fixed point results. An a priori estimate in $L^\infty(\R^N) \times L^\infty(\R^N)$ for solutions of \eqref{probwholeone} is first established by a Moser's type iteration procedure and an adequate truncation, which, due to singular terms, require a specific treatment. Problem \eqref{probwholeone} is next perturbed by introducing a parameter $\epsilon>0$. This produces the family of regularized systems
\begin{equation}\label{pr}
\left\{
\begin{alignedat}{2}
-\Delta_p u =a_1(x)f(u+\epsilon,v)\;\; & \mbox{in $\R^N$,}\;\; u>0\;\; && \mbox{in $\R^N$,}\\
-\Delta_q v =a_2(x)g(u,v+\epsilon)\;\; & \mbox{in $\R^N$,}\;\; v>0\;\; && \mbox{in $\R^N$,}
\end{alignedat}
\right.
\end{equation}
whose study yields useful information on the  original problem. In fact, the previous $L^\infty$-boundedness still holds for solutions to \eqref{pr}, regardless of $\epsilon$. Thus, via Schauder's fixed point theorem, one gets a solution $(u_\epsilon,v_\epsilon)$ lying inside a rectangle given by positive lower bounds, where $\epsilon$ does not appear, and positive upper bounds, that may instead depend on $\epsilon$. Finally, letting $\epsilon\to 0^+$ and using the $({\rm S})_+$-property of the negative $r$-Laplacian in $\mathcal{D}^{1,r}(\R^N)$ (see \cite[Proposition 2.2]{MMM}) yields a weak solution to \eqref{probwholeone}.
\begin{thm}[\cite{MMM}, Theorem 5.1]
Let $({\rm a}_{12})$ and $({\rm a}_{13})$ be satisfied. Then \eqref{probwholeone} has a weak solution $(u,v)\in X^{p,q}(\R^N)$, which is essentially bounded.
\end{thm} 
Very recently, the parametric system
\begin{equation}\label{probwholetwo}
\left\{
\begin{alignedat}{2}
-\Delta_p u =a_1(x)f_1(u)+\lambda b_1(x)g_1(u)h_1(v)\;\; & \mbox{in $\R^N$,}\;\; u>0\;\; && \mbox{in $\R^N$,}\\
-\Delta_p v =a_2(x)f_2(v)+\mu b_2(x)g_2(v)h_2(u)\;\; & \mbox{in $\R^N$,}\;\; v>0\;\; && \mbox{in $\R^N$,}\\
u(x)\to 0,\; v(x)\to 0\;\;\mbox{as $|x|\to\infty$}
\end{alignedat}
\right.
\end{equation}
where $N\geq 3$, $1<p<N$, $a_i,b_i\in C^0(\R^N)$, $f_i,g_i,h_i\in C^0(\R^+,\R^+)$, $f_i$ is singular at zero, and $\lambda,\mu>0$, was thoroughly investigated in \cite{SLRZ}. Under suitable hypotheses, it is shown that there exists an open set $\Theta\subseteq (\R^+)^2$, whose part of its boundary contained in $(\R^+)^2$, say $\Gamma$, turns our a continuous monotone curve, such that \eqref{probwholetwo} admits a $C^1$-solution if $(\lambda,\mu)\in\Theta$ and has no solution when $(\lambda,\mu)\in(\R^+)^2\setminus (\Theta\cup\Gamma)$.
\subsection{Uniqueness}
As far as we know, uniqueness of solutions to singular quasi-linear elliptic systems in the whole space is still an open problem. Taking inspiration from \cite{CDS}, a first result has been obtained by Gambera and Guarnotta \cite{GG1}.
\subsection{Systems with convection terms}
The very recent paper \cite{GMM} treats the problem
\begin{equation}\label{probwholethree}
\left\{
\begin{alignedat}{2}
-\Delta_p u=f(x,u,v,\nabla u,\nabla v)\;\; &\mbox{in $\R^N$,}\;\; u>0\;\; &&\mbox{in $\R^N$,}\\
-\Delta_q v=g(x,u,v,\nabla u,\nabla v)\;\; &\mbox{in $\R^N$,}\;\; v>0\;\; &&\mbox{in $\R^N$,}
\end{alignedat}
\right.
\end{equation}
where $N\geq 3 $, $p,q\in ]2-\frac{1}{N},N[$, while $f,g:\R^N \times(\R^+)^2\times \R^{2N} \to\R^+$ are Carathéodory's functions satisfying assumptions $({\rm a}_{14})$--$({\rm a}_{15})$ below. Since $f,g$ depend on the gradient of solutions and equations are set in the whole space, neither variational methods can be exploited nor compactness for Sobolev embedding holds. The research started in \cite{MMM}, where convective terms did not appear, is continued here, along the works \cite{CLM,GM,GMMotr}, which address analogous questions, but concerning a bounded domain. 

A pair $(u,v)\in X^{p,q}(\R^N)$ such that $u,v>0$  is called:
\begin{itemize}
\item[1)] \emph{distributional solution} to \eqref{probwholethree} if for every $(\phi,\psi)\in C^{\infty}_0(\R^N)^2$ one has
\begin{equation}\label{intbyparts}
\begin{split}
\int_{\R^N}|\nabla u|^{p-2}\nabla u\nabla\phi\,{\rm d}x &= \int_{\R^N} f(\cdot,u,v,\nabla u,\nabla v) \phi\,{\rm d}x, \\
\int_{\R^N}|\nabla v|^{q-2}\nabla v\nabla\psi\,{\rm d}x &= \int_{\R^N} g(\cdot,u,v,\nabla u,\nabla v) \psi\,{\rm d}x;
\end{split}
\end{equation}
\item[2)] \emph{(weak) solution} of \eqref{probwholethree} when \eqref{intbyparts} holds for all $(\phi,\psi)\in X^{p,q}(\R^N)$;
\item[3)] \emph{`strong' solution} to \eqref{probwholethree} if $|\nabla u|^{p-2}\nabla u, |\nabla v|^{q-2}\nabla v \in W^{1,2}_{\rm loc}(\R^N)$ and the differential equations are satisfied a.e. in $\R^N$.
\end{itemize}
Obviously, both 2) and 3) force 1), whilst reverse implications turn out generally false; see also \cite[Remark 4.5]{GMM}. Moreover, as observed at p. 48 of \cite{SS}, problems in unbounded domains may admit strong solutions that are not weak or vice-versa.

Roughly speaking, the technical approach proceeds as follows. The auxiliary problem
\begin{equation}\label{probwholeaux}
\left\{
\begin{alignedat}{2}
-\Delta_p u=f(x,u+\epsilon,v,\nabla u,\nabla v)\;\; &\mbox{in $\R^N$,}\;\; u>0\;\; &&\mbox{in $\R^N$,}\\
-\Delta_q v=g(x,u,v+\epsilon,\nabla u,\nabla v)\;\; &\mbox{in $\R^N$,}\;\; v>0\;\; &&\mbox{in $\R^N$,}
\end{alignedat}
\right.
\end{equation}
$\epsilon>0$, obtained by shifting appropriate variables of reactions, which avoids singularities, is first solved. To do this, nonlinear regularity theory, a priori estimates, Moser's iteration, trapping region, and fixed point arguments are employed. Unfortunately, bounds from above alone do not allow to get a solution of \eqref{probwholethree}: treating singular terms additionally requires some estimates from below. Theorem 3.1 in \cite{DaMi} ensures that solutions to \eqref{probwholeaux} turn out locally greater than a positive constant regardless of $\epsilon$. Thus, under the hypotheses below, one can construct a sequence $\{(u_\epsilon,v_\epsilon)\}\subseteq X^{p,q}(\R^N)$ such that $(u_\epsilon,v_\epsilon)$ solves \eqref{probwholeaux} for all $\epsilon>0$ and whose weak limit as $\epsilon\to 0^+$ is a distributional solution to \eqref{probwholethree}. Next, a localization-regularization reasoning shows that 
$$(u,v)\;\mbox{distributional solution}\implies (u,v)\;\mbox{weak solution.}$$
Through a recent differentiability result \cite[Theorem 2.1]{CM} one then has
$$(u,v)\;\mbox{distributional solution}\implies (u,v)\;\mbox{strong solution.}$$
The assumptions below will be posited. 
\begin{itemize}
\item[$({\rm a}_{14})$] There exist $\alpha_i\in ]-1,0]$, $\beta_i,\gamma_i, \delta_i\in\R^+_0$, as well as $\mu_i,\hat{\mu}_i> 0$  such that
\begin{equation*}
\begin{aligned}
\mu_1 a_1(x) s^{\alpha_1} t^{\beta_1}\leq f(x,s,t,\xi_1,\xi_2)\leq\hat{\mu}_1 a_1(x)
\left( s^{\alpha_1} t^{\beta_1} + |\xi_1|^{\gamma_1} + |\xi_2|^{\delta_1} \right)\\
\mu_2 a_2(x) s^{\beta_2} t^{\alpha_2}\leq g(x,s,t,\xi_1,\xi_2)\leq\hat{\mu}_2 a_2(x)
\left( s^{\beta_2} t^{\alpha_2} + |\xi_1|^{\gamma_2} + |\xi_2|^{\delta_2} \right)
\end{aligned}
\end{equation*}
in $\R^N\times(\R^+)^2\times\R^{2N}$. Moreover,
$$\beta_1\vee\delta_i<q-1,\;\;\;\beta_2\vee \gamma_i<p-1,\;\; i=1,2,$$
$a_1\in L^{s_p}_{\rm loc}(\R^N)$, with $s_p >p'N$, $a_2\in L^{s_q}_{\rm loc}(\R^N)$, with $s_q >q'N$, and $\displaystyle {\essinf_{B_\rho} a_i > 0}$ for all $\rho>0$.

\item[$({\rm a}_{15})$] There exist $\zeta_1,\zeta_2\in]N,+\infty] $ such that 
$a_i\in L^1(\R^N)\cap L^{\zeta_i}(\R^N)$, where
\begin{equation*}
\frac{1}{\zeta_1} < 1 - \frac{p}{p^*} - \theta_1, \quad \frac{1}{\zeta_2} < 1 - \frac{q}{q^*} - \theta_2,
\end{equation*}
with
\begin{equation*}
\theta_1:=\max\left\{\frac{\beta_1}{q^*},\frac{\gamma_1}{p},\frac{\delta_1}{q}\right\}< 1-\frac{p}{p^*},\quad
\theta_2:=\max\left\{\frac{\beta_2}{p^*},\frac{\gamma_2}{p},\frac{\delta_2}{q}\right\}< 1-\frac{q}{q^*}.
\end{equation*}
Further,
\begin{equation*}
\begin{aligned}
(\beta_1\vee\delta_1)(\beta_2\vee\gamma_2)< (p-1-\gamma_1)(q-1-\delta_2),\\
\phantom{}\\
\frac{1}{s_p}+\left(\frac{\gamma_1}{p}\vee\frac{\delta_1}{q}\right)\leq \frac{1}{2}, \qquad
\frac{1}{s_q}+\left(\frac{\gamma_2}{p}\vee\frac{\delta_2}{q}\right)\leq \frac{1}{2}.
\end{aligned}
\end{equation*}
\end{itemize}
\begin{ex}
Condition $({\rm a}_{15})$ is fulfilled once $a_1,a_2 \in L^1(\R^N)\cap L^\infty(\R^N)$ and
\begin{equation*}
\max \left\{ \frac{\beta_1}{q^*},\frac{\gamma_1}{p},\frac{\delta_1}{q} \right\}< 1-\frac{p}{p^*}, \quad
\max \left\{ \frac{\beta_2}{p^*},\frac{\gamma_2}{p},\frac{\delta_2}{q} \right\}< 1-\frac{q}{q^*}.
\end{equation*}
In fact, it suffices to choose $ \zeta_1:= \zeta_2:= +\infty $.
\end{ex}
\begin{thm}[\cite{GMM}, Theorem 1.3]
Under hypotheses $({\rm a}_{14})$--$({\rm a}_{15})$, problem \eqref{probwholethree} admits a weak and strong solution $(u,v)\in X^{p,q}(\R^N)$.
\end{thm}
\begin{rmk}
If we merely seek weak solutions to \eqref{probwholethree} then the request $p,q\in]1,N[$ and a weaker integrability property of $a_i$ suffice; cf. \cite[Section 4.2.2]{Gdr}.
\end{rmk}
\small{
\subsection*{Acknowledgement}
U. Guarnotta and S.A. Marano were supported by the research project PRA 2020--2022, Linea 3 and Linea 2, resp., `MO.S.A.I.C.' of the University of Catania.\\
U. Guarnotta, R. Livrea, and S.A. Marano were supported by the research project PRIN 2017 `Nonlinear Differential Problems via Variational, Topological and Set-valued Methods' (Grant No. 2017AYM8XW) of MIUR.
}


\begin{thebibliography}{777}
%
\bibitem{AC}
C.O. Alves and F.J.S.A. Correa, \textit{On the existence of positive solution for a class of singular systems involving quasilinear operators,} Appl. Math. Comput. \textbf{185} (2007), 727--736.
%
\bibitem{AM}
C.O. Alves and A. Moussaoui, \textit{Existence of solutions for a class of singular elliptic systems with convection term,} Asymptot. Anal. \textbf{90} (2014), 237--248.
%
\bibitem{AM1}
C.O. Alves and A. Moussaoui \textit{Existence and regularity of solutions for a class of singular $(p(x),q(x))$-Laplacian systems,} Complex Var. Elliptic Equ. (2017), DOI 10.1080/17476933.2017.1298589.
%
\bibitem{AMT}
C.O. Alves, A. Moussaoui, and L. Tavares, \textit{An elliptic system with logarithmic
nonlinearity,}Adv. Nonlinear Anal. \textbf{8} (2019), 928--845. 
%
\bibitem{AN}
S.C. Arruda and R.G. Nascimento, \textit{Existence and multiplicity of positive solutions for a singular system via sub-supersolution method and Mountain Pass Theorem,} Electron. J. Qual. Theory Differ. Equ. 2021, Paper no. 26, 20 pp.
%
\bibitem{B1}
M. Benrhouma, \textit{Existence and uniqueness of solutions for a singular semilinear system,} Nonlinear Anal. \textbf{107} (2014), 134--146.
%
\bibitem{B2}
M. Benrhouma, \textit{On a singular elliptic system with quadratic growth in the gradient,} J. Math. Anal. Appl. \textbf{448} (2017), 1120--1146.
%
\bibitem{BO}
L. Boccardo and L. Orsina, \textit{A variational semilinear singular system,} Nonlinear Anal. \textbf{74} (2011), 3849--3860.
%
\bibitem{CLM}
P. Candito, R. Livrea, and A. Moussaoui, \textit{Singular quasilinear elliptic systems involving gradient terms,} Nonlinear Anal. Real World Appl. \textbf{55} (2020), 103142.
%
\bibitem{CaMo}
S. Carl and D. Motreanu. \textit{Extremal solutions for novariational quasilinear elliptic systems via expanding trapping regions,} Monatsh. Math. \textbf{182} (2017), 801--821.
%
\bibitem{CMS}
J. Carmona, P.J. Martinez-Aparicio, and A. Suarez, \textit{A sub-supersolution method for nonlinear elliptic singular systems with natural growth and some applications,} Nonlinear Anal. \textbf{132} (2016), 47--65.
%
\bibitem{CSGS}
M.L.M. Carvalho, E.D. Silva, C. Goulart, and C.A. Santos, 
\textit{Ground and bound state solutions for quasilinear elliptic systems including singular nonlinearities and indefinite potentials,} Comm. Pure Appl. Math. \textbf{19} (2020), 4401--4432.
%
\bibitem{CDS}
M. Chhetri, P. drabek, and R. Shivaji, \textit{Analysis of positive solutions for classes of quasilinear singular problems on exterior domains,} Adv. Nonlinear Anal. \textbf{6} (2017), 447--459.
%
\bibitem{CMcK1}
Y.S. Choi and P.J. McKenna, \textit{A singular Gierer-Meinhardt system of elliptic equations,} Ann. Inst. H. Poincar\'{e} Anal. Non Lin\'{e}aire \textbf{17} (2000), 503--522.
%
\bibitem{CMcK2}
Y.S. Choi and P.J. McKenna, \textit{A singular Gierer-Meinhardt system of elliptic equations: the classical case,} Nonlinear Anal. \textbf{55} (2003), 521--541.
\bibitem{CHS}
K.D. Chu, D.D. Hai, and R. Shivaji, \textit{Positive solutions for a class of non-cooperative $pq$-Laplacian systems with singularities,} Appl. Math. Lett. \textbf{85} (2018), 103--109.
%
\bibitem{CM}
A. Cianchi and V.G. Maz'ya, \textit{Second-order two-sided estimates in nonlinear elliptic problems,} Arch. Ration. Mech. Anal. \textbf{229} (2018), 569–599.
%
\bibitem{Calt}
F.J.S.A. Correa, G.C.G. dos Santos, L.S. Tavares, and S.S, Muhassua, \textit{Existence of solution for a singular elliptic system with convection terms,} Nonlinear Aanal. Real World Appl. \textbf{66} (2022), Paper no. 103549, 18 pp.
%
\bibitem{DaMi}
L. D'Ambrosio and E. Mitidieri, \textit{Entire solutions of quasilinear elliptic systems on Carnot groups,} Reprint of Tr. Mat. Inst. Steklova \textbf{283} (2013), 9–24, Proc. Steklov Inst. Math. \textbf{283} (2013), 3–19.
%
\bibitem{DeC}
L. M. De Cave, \textit{Singular elliptic systems with higher order terms of p-laplacian type,} Adv. Nonlinear Stud. \textbf{16} (2016), 667--687.
%
\bibitem{DOS}
L.M. De Cave, F. Oliva, and M. Strani, \textit{Existence of solutions to a non-variational singular elliptic system with unbounded weights,} Math. Nachr. \textbf{290} (2017), 236–247.
%
\bibitem{DeMou}
H. Dellouche and A. Moussaoui, \textit{Singular quasilinear elliptic systems with gradient dependence,} Positivity \textbf{26} (2022), paper no. 10, 14 pp.

\bibitem{DKC} 
M. del Pino, M. Kowalczyk, and X. Chen, \emph{The Gierer-Meinhardt system: the breaking of homoclinics and multi-bump ground states,} Commun. Contemp. Math. \textbf{3} (2001), 419--439.
%
\bibitem{DKW} 
M. del Pino, M. Kowalczyk, and J. Wei, \emph{Multi-bump ground states of the Gierer-Meinhardt system in $\R^2$,} Ann. Inst. H. Poincar\'{e}, Anal. Non Lin\'{e}aire \textbf{20} (2003), 53--85.
%
\bibitem{DL}
J. Deny and J.L. Lions, \textit{Les espaces du type de Beppo Levi,} Ann. Inst. Fourier (Grenoble) \textbf{5} (1955), 305–370.
%
\bibitem{DM}
H. Didi and A. Moussaoui, \textit{Multiple positive solutions for a class of quasilinear singular elliptic systems,} Rend. Circ. Mat. Palermo (2) \textbf{69} (2020), 977--994.
%
%
\bibitem{EPS}
S. El Manouni, K. Perera, and R. Shivaji; \textit{On a singular quasi-monotone $(p,q)$-Laplacian systems,} Proc. Roy. Soc. Edinburgh Sect. A \textbf{142} (2012), 585--594. 
%
\bibitem{Ga}
G.P. Galdi, \textit{An introduction to the mathematical theory of the Navier-Stokes equations. Steady-state problems,} 2nd ed., Springer Monographs in Mathematics, Springer, New York, 2011.
%
\bibitem{GG1}
L. Gambera and U. Guarnotta, in preparation.
%
\bibitem{Ghe1}
M. Ghergu, \textit{Singular elliptic systems of Lane-Emden type,} in ``Recent trends in nonlinear partial differential equations. II. Stationary problems'', 253--262, Contemp. Math., \textbf{595}, Amer Math. Soc., Providence, RI, 2013.
%
\bibitem{Ghe}
M. Ghergu, \textit{On a class of singular elliptic systems,} Nonlinear Anal. \textbf{119} (2015), 98--105.
%
\bibitem{GR}
M. Ghergu and V. Radulescu, \textit{On a class of singular Gierer-Meinhardt systems arising in morphogenesis,} C. R. Math. Acad. Sci. Paris \textbf{344} (2007), 163--168.
%
\bibitem{GR1}
M. Ghergu and V.D. Radulescu, \textit{Singular elliptic problems: bifurcation and asymptotic analysis}, Oxford Lecture Ser. Math. Appl. \textbf{37}, Oxford Univ. Press, Oxford, 2008.
%
\bibitem{GHM}
J. Giacomoni, J. Hernandez, and A. Moussaoui, \textit{Quasilinear and singular systems: the cooperative case,} Contemp. Math. \textbf{540} (2011), 79--94.
%
\bibitem{GHS}
J. Giacomoni, J. Hernandez, and P. Sauvy, \textit{Quasilinear and singular elliptic systems,} Adv. Nonlinear Anal. \textbf{2} (2013), 1--41.
%
\bibitem{GST}
J. Giacomoni, I. Schindler, and P. Takac, \textit{Singular quasilinear elliptic systems and H\"{o}lder regularity,} Adv. Differential Equations \textbf{20} (2015), 259--298.
%
\bibitem{GiMe}
A. Gierer and H. Meinhardt, \textit{Generation and regeneration of sequence of structures during morphogenesis,} J, Theor. Biol. \textbf{85} (1980), 429--450.
%
\bibitem{Go}
T. Godoy, \textit{Existence of positive weak solutions for a nonlocal singular elliptic system,} AIMS Math. \textbf{4} (2019), 792--804.
%
\bibitem{Gdr}
U. Guarnotta, \textit{Existence results for singular convective elliptic problems} (PhD thesis), Department of Mathematics and Computer Sciences, University of Palermo (Italy), AY 2020/21.
%
\bibitem{GM}
U. Guarnotta and S.A. Marano, \textit{Infinitely many solutions to singular convective Neumann systems with arbitrarily growing reactions,} J. Differential Equations \textbf{271} (2021), 849–863.
%
\bibitem{GM1}
U. Guarnotta and S.A. Marano, \textit{Corrigendum to “Infinitely many solutions to singular convective Neumann systems with arbitrarily growing reactions” [J. Differ. Equ. 271 (2021) 849–863],} J. Differential Equations \textbf{274} (2021), 1209--1213.
%
\bibitem{GMMotr}
U. Guarnotta, S.A. Marano, and D. Motreanu, \textit{On a singular Robin problem with convection terms,} Adv. Nonlinear Stud. \textbf{20} (2020), 895–909.
%
\bibitem{GMM} 
U. Guarnotta, S.A. Marano, and A. Moussaoui,
\textit{Singular quasilinear convective elliptic systems in $\R^N$,} Adv. Nonlinear Anal. \textbf{11} (2022), 741–756.
%
\bibitem{HS}
D.D. Hai and R.C. Smith, \textit{Uniqueness for a class of singular semilinear elliptic systems,} Funkcial. Ekvac. \textbf{59} (2016), 35--49.
%
\bibitem{HMV}
J. Hern\'{a}ndez, F.J. Mancebo, and J.M. Vega, \textit{Positive solutions for singular semilinear elliptic systems,} Adv. Differential Equations \textbf{13} (2008), 857--880.
%
\bibitem{Ka}
O. Kavian, \textit{Inegalit\'{e} de Hardy–Sobolev et applications,} Th\`{e}se de Doctorate de 3eme cycle, Universit\'{e} de Paris VI (1978).
%
\bibitem{KM}
B. Khodja and A. Moussaoui, \textit{Positive solutions for infinite semipositone/positone quasilinear elliptic systems with singular and superlinear terms,}
Differ. Equ. Appl. \textbf{8} (2016), 535--546.
%
\bibitem{Kim}
E.H. Kim, \textit{Singular Gierer-Meinhardt systems of elliptic boundary value problems,} J. Math. Anal. Appl. \textbf{308} (2005), 1--10.
%
\bibitem{K}
M.A. Krasnoselskii, \textit{Topological Methods in the Theory of Nonlinear Integral Equations,} Pergamon Press, Oxford-London-Paris, 1964.
%
\bibitem{KS}
M. Kratou and K. Saoudi, \textit{The fibering map approach for a singular elliptic system involving the $p(x)$-laplacian and nonlinear boundary conditions.} Rev. Un. Mat. Argentina \textbf{62} (2021), 171--189.
%
\bibitem{LSY}
E.K. Lee, R. Shivaji, and J. Ye, \textit{Classes of singular $pq$-Laplacian semipositone systems,} Discrete Contin. Dyn. Syst. \textbf{27} (2010), 1123--1132. 
%
\bibitem{LL}
E.H. Lieb and M. Loss, \textit{Analysis,} Graduate Studies in Mathematics, vol. 14, 2nd ed., American Mathematical Society, Providence, 2001.
%
\bibitem{L}
G.M. Lieberman, \textit{Boundary regularity for solutions of degenerate elliptic equations,} Nonlinear Anal. \textbf{12} (1988), 1203–1219.
%
\bibitem{MMM}
S.A. Marano, G. Marino, and A. Moussaoui, \textit{Singular quasilinear elliptic systems in $\R^N$,} Ann. Mat. Pura Appl. \textbf{198} (2019), 1581--1594.
%
\bibitem{Mi}
L. Mi, \textit{Positive solutions for a class of singular elliptic systems,} Electron. J. Qual. Theory Differ. Equ. \textbf{24} (2017), 1--13.
%
\bibitem{MS}
M. Montenegro and A. Su\'arez, \textit{Existence of a positive solution for a singular system,} Proc. Roy. Soc. Edinburgh Sect. A \textbf{140} (2010), 435--447.
%
\bibitem{Mo}
D. Motreanu, \textit{Nonlinear differential problems with smooth and nonsmooth constraints,} Math. Anal. Appl. Ser., Academic Press, London, 2018.
%
\bibitem{MM1}
D. Motreanu and A. Moussaoui, \textit{Existence and boundedness of solutions for a singular cooperative quasilinear elliptic system,} Complex Var. Elliptic Equ. \textbf{59} (2014), 285--296.
%
\bibitem{MM2}
D. Motreanu and A. Moussaoui, \textit{An existence result for a class of quasilinear singular competitive elliptic systems,} Appl. Math. Lett.\textbf{38} (2014), 33--37.
%
\bibitem{MM3}
D. Motreanu and A. Moussaoui, \textit{A quasilinear singular elliptic system without cooperative structure,} Acta Math. Sci. Ser. B Engl. Ed. \textbf{34} (2014), 905--916.
%
\bibitem{MMZ}
D. Motreanu, A. Moussaoui, and Z. Zhang, \textit{Positive solutions for singular elliptic systems with convection term,} J. Fixed Point Theory Appl. \textbf{19} (2017), 2165--2175.
%
\bibitem{MKT}
A. Moussaoui, B. Khodja, and S. Tas, \emph{A singular Gierer-Meinhardt system of elliptic equations in $\R^N$,} Nonlinear Anal. \textbf{71} (2009), 708--716.
%
\bibitem{MV}
A. Mossaoui and J. Velin, \textit{Existence and a priori estimates of solutions for quasilinear singular elliptic systems with variable exponents,} J. Elliptic Parabol. Equ. \textbf{4} (2018), 417--440.
%
\bibitem{PuSe}
P. Pucci and J. Serrin, \textit{The maximum principle,} Prog. Nonlinear Differential Equations Appl. {\bf 73}, Birkh\"auser Verlag, Basel, 2007.
%
\bibitem{Ra}
P.H. Rabinowitz, \textit{Some global results for nonlinear eigenvalue problems,} J. Funct. Anal. {\bf 7} (1971) 487–513.
%
\bibitem{SLRZ}
C.A. Santos, R. Lima Alves, M. reis, and J. Zhou, \textit{Maximal domains of the
$(\lambda,\mu)$-parameters to existence of entire positive solutions for
singular quasilinear elliptic systems,} J. Fixed Point Theory Appl. \textbf{22} (2020), Paper no. 54, 30 pp.
%
\bibitem{S}
K. Saoudi, \textit{A singular elliptic system involving the $p(x)$-Laplacian and generalized Lebesgue–Sobolev spaces,} Internat. J. Math. (2019), Paper no. 1950064, 17 pp.
%
\bibitem{SS}
C.G. Simader and H. Sohr, \textit{The Dirichlet problem for the Laplacian in bounded and unbounded domains. A new approach to weak, strong and $ (2+k) $-solutions in Sobolev-type spaces,} Pitman Research Notes in Mathematics Series \textbf{360}, Longman, Harlow, 1996.
%
\bibitem{Si}
G. Singh, \textit{Weak solutions for singular quasilinear elliptic systems,} Complex Var. Elliptic Equ. \textbf{61} (2016), 1389--1408.
%
\end{thebibliography}
\end{document}